\documentclass[11pt, oneside, leqno]{amsart}
\usepackage{graphicx}
\usepackage{amssymb, amsmath}

\input xy
\xyoption{all}
\CompileMatrices

\newcommand{\bo}[1]{\mathbf{#1}}

\newcommand{\ra}{\rightarrow}
\newcommand{\xra}{\xrightarrow}

\newcommand{\gr}{\bo {Groups}}
\newcommand{\Hom}{{\rm Hom}}

\newcommand{\ralim}{\varinjlim}

\newcommand{\Q}{{\mathbb{Q}}}
\newcommand{\Z}{{\mathbb{Z}}}

\newtheorem{theorem}{Theorem}[section]
\newtheorem{lemma}[theorem]{Lemma}

\newtheorem{proposition}[theorem]{Proposition}

\newenvironment{Proof}{\removelastskip\par\medskip
\noindent{\em Proof.} \rm}{\penalty-20\null\hfill${\,\lower0.9pt
\vbox{\hrule \hbox{\vrule height 0.3 cm \hskip 0.3 cm \vrule height 0.3 cm}
\hrule}\,} $\par\medbreak}

\swapnumbers

\title{ A note on localizations of perfect groups}

\author[B. Badzioch]{Bernard Badzioch}
\address{
University of Minnesota\\
School of Mathematics\\
Minneapolis, MN 55455\\
}
\email[B. Badzioch]{badzioch@math.umn.edu}
\author[M. Feshbach]{Mark Feshbach}
\email[M. Feshbach]{feshbach@math.umn.edu}

\date{January 10, 2003}

\begin{document}

\begin{abstract}
We describe a perfect group whose localization is not perfect. 
\end{abstract}

\maketitle

\section{Introduction}
A localization is a type of a functor $L\colon \gr\ra\gr$ which is idempotent 
(that is  $LLG\cong LG$) and admits  a coaugmentation $\eta\colon G\ra LG$ \cite{MR2001i:55012}.   
Localizations are ubiquitous in group theory: abelianization, killing of the $p$--torsion and 
inversion of a prime in a group are all examples of  such functors.  

The natural question --  which classes of of groups are preserved by all localizations -- has been 
 a focus of a lot of study recently. This work yield both classes which are preserved
(e.g. abelian groups, nilpotent groups of class 2  \cite{MR2001c:20107}), 
and these which do not have this property, 
like finite \cite{MR2001c:20108} or solvable groups \cite{MR2001k:20070}. 
The goal of this note is to prove the following 
\begin{theorem} \label{MAIN}
The class of perfect groups is not closed with respect to taking localizations. 
That is,  there exists a perfect group $P$ and a localization $\eta\colon P\rightarrow LP$ 
such that $LP$ is not perfect. 
\end{theorem}
\noindent This answers a question posed by Casacuberta in \cite{MR2001i:55012}.   

The groups  $P$ and $LP$ we construct in the proof of Theorem \ref{MAIN} are infinite. 
It would be interesting to know if one can find a finite perfect group with a non-perfect localization. 
The following shows however that the localized group would have to be infinite. 
\begin{proposition} \label{FINITE}
If $\eta\colon P\ra LP$ is a localization of a perfect group $P$ and $LP$ is finite then
$LP$ is a perfect group. 
\end{proposition}


\section{Proofs}
Our main tool will be the following fact which characterizes all possible localizations of
a group. It is a direct ramification  of the definition of  localization functors 
(see e.g. \cite[Lemma 2.1]{MR2001i:55012}).

\begin{lemma} \label{HOM LOC}
A homomorphism $\eta\colon G\ra  H$ is a localization of $G$ with respect to some localization 
functor iff  $\eta$ induces a bijection of sets 
$$\Hom(H, H)\xra{\eta^\ast}\Hom(G ,H)$$
\end{lemma}  

\noindent As an application we obtain 
\begin{lemma}\label{HOM PERF}
If $\eta\colon G\ra H$ is a localization and $G$ is a perfect group then there are
no non-trivial homomorphisms $H/[H,H]\ra H$.   
\end{lemma}

\begin{proof}
Let  $g\colon H/[H,H]\ra H$ be any homomorphism. 
and let $f$ denote the composition $H\ra H/[H, H]\xra{g} H$. 
Since $G$ is perfect the composition $f\circ \eta$ is the trivial map. 
Lemma \ref{HOM LOC} implies then that  $f$ is also trivial, and thus so is $g$.     
\end{proof}

\noindent Since every finite group $H$ admits a non-trivial  map $H/[H,H]\ra H$
unless $H$ is perfect,  Proposition \ref{FINITE} is a consequence of Lemma \ref{HOM LOC}.  

Next, we turn to the proof of Theorem \ref{MAIN}.  We start with


\vskip .5cm

\noindent{\bf Construction of the group P. }
For $n\geq 0$ let $\tilde P_n$ be a free group on 
$2^n$ generators $x^{(n)}_1$, $x^{(n)}_2$,\dots,$x^{(n)}_{2^n}$, and 
let $\tilde\phi_n\colon \tilde P_n\rightarrow \tilde P_{n+1}$
be a group homomorphism defined by 
$$\tilde \phi_n(x^{(n)}_i)=[x^{(n+1)}_{2i-1}, x^{(n+1)}_{2i}]$$
where $[a,b]=a^{-1}b^{-1}ab$ is the commutator of $a$ and $b$. 
Define $\tilde P := \ralim_{n} \tilde P_n$. Notice,  that since 
$\tilde P$ is generated by the elements 
$x^{(n)}_i\in [\tilde P, \tilde P]$ the group $\tilde P$ is perfect. 
Let $K$ be the smallest normal subgroup of $\tilde P$ containing the 
elements $[x^{(0)}_1, x^{(n)}_i]$ for all $n\geq 0$, $1\leq i \leq 2^n$. 
Define $P:= \tilde P/K$.

\begin{proposition}
\label{P PERFECT}
The group $P$ is perfect and $x^{(0)}_1$ is a central element of P. 
Moreover, $x^{(0)}_1$ is an element of infinite order,  and as a consequence 
$P$ is a non-trivial group. 
\end{proposition}

\begin{proof}
The first two claims are obvious. To see that $x^{(0)}_1\in P$ 
has infinite order notice that P can be viewed as a limit 
$$P=\ralim_n P_n$$
where $P_n$ is a group with the presentation
$$P_n:=\langle x^{(n)}_1, \dots, x^{(n)}_{2^n} \   | \  [x^{(0)}_1,x^{(n)}_{i}]=1,
\ \ i=1,\dots, 2^n \rangle$$
(by abuse of notation we denote here by $x^{(0)}_1$ the image of 
the element $x^{(0)}_1$ under the map $P_0\rightarrow P_n$). 
It is then enough to show that $x^{(0)}_1$ has infinite order in 
$P_n$ for all $n\geq 0$. To see this consider  $GL(\Z, 2^n+1)$  -- the group 
of invertible matrices of dimension 
$2^{n}+1$ with integer coefficients. For $n>0$ there is a homomorphism
$$\psi_n\colon P_n\rightarrow  GL(\Z, 2^n+1)$$
defined by $\psi_n(x^{(n)}_i)= e^{1}_{i,i+1}$, where $e^a_{i,j}$ denote 
the matrix with 1's on the diagonal, $a$ as the $(i,j)$-th entry, and 
and 0's elsewhere. One can check that $\psi_n(x^{(0)}_1)=e^{\pm 1}_{1,2^n+1}$. 
Since $(e^{\pm 1}_{1,2^n+1})^k=e^{\pm k}_{1,2^n+1}$ this is a non-torsion element 
of $GL(\Z, 2^n+1)$,  and  as a consequence  $x^{(0)}_1$ has infinite order 
in  $P_n$ as claimed. 
\end{proof}


\vskip .5cm 

\noindent {\bf Construction of the map $\bf \boldsymbol{ \eta}\colon P\boldsymbol{\ra }LP$. }
Let $\Q$ be the group of rational numbers. Define 
$$LP:=P\oplus \Q/\langle(x^{(0)}_1, -1)\rangle$$
and let the map $\eta\colon P\rightarrow LP$ be given by the composition 
of the inclusion $P\hookrightarrow P\oplus \Q$ and the projection 
$P\oplus \Q\rightarrow LP$. Since $\eta$ is a monomorphism we will 
identify $P$ with its image $\eta(P)$. Notice that $P$ is a normal subgroup 
of $LP$ and that $LP/P\cong \Q/\Z$. Since $\Q/\Z$ is not a perfect group 
neither is $LP$. 

It remains to prove that $\eta$ is a localization of $P$. By Lemma \ref{HOM LOC} 
this amount showing that any homomorphism $f\colon P\rightarrow LP$ admits 
a unique factorization
$$\xymatrix{
P\ar[r]^{\eta}\ar[dr]_f & LP\ar[d]^{\bar f}\\
& LP\\
}$$

\vskip .5cm 

\noindent {\bf Uniqueness of $\bf \bar f$. }
Assume that $\bar f_1, \bar f_2 \colon LP\rightarrow LP$
are homomorphisms such that $\eta\bar f_1 = \eta \bar f_2$, and consider 
the homomorphism 
$$g:= (\bar f_1|_\Q - \bar f_2|_\Q)\colon \Q \rightarrow LP$$
We have $\bar f_1(1)=\bar f_1(x^{(0)}_1)=\bar f_2(x^{(0)}_1)=\bar f_2(1)$,
and thus $\Z\subseteq \ker g$. Therefore we get a factorization 
$$g\colon \Q/\Z\rightarrow LP$$
and $g\equiv 1$ iff $\bar f_1=\bar f_2$. Thus, our claim is a consequence of
the following

\begin{lemma}
\label{LP TORSION}
The group $LP$ is torsion free. 
\end{lemma}

\begin{proof}
Let $(w, \frac{p}{q})$ represents a torsion element of $LP$.
Then $q\cdot(w, \frac{p}{q})=(w^q, p)= w^q(x^{(0)}_1)^p$ is a torsion 
element in P. Consider the group $R:=P/\langle x^{(0)}_1\rangle$. The element 
$ w^q(x^{(0)}_1)^p=w^q$ is  torsion in $R$, and thus so is $w$. 
Notice that $R=\ralim_n R_n$ where 
$$R_n=\langle x^{(n)}_1,\dots, x^{(n)}_{2^n} | x^{(0)}_1 =1\rangle$$
It follows that $w$ is a torsion element in $R_n$ for  $n$ large enough. 
On the other hand, $R_n$ is a group with one relator given by a word 
which is not a proper power of any element in the free group. By    
\cite[Thm. 4.12, p. 266]{MR34:7617}  $R_n$ must be  torsion free. 
Therefore 
$w=1$ in R, and  so $(w, \frac{p}{q})= ((x^{(0)}_1)^l, \frac{p}{q})= 
l+\frac{p}{q}\in \Q\subseteq LP$ for some $l\in \Z$. Since by 
assumption  $(w, \frac{p}{q})$ is a torsion element it must be trivial. 
 
\end{proof}


\vskip .5cm 

\noindent {\bf Existence of $\bf \bar f$.}
We need to show that  every  homomorphism $f\colon P\rightarrow LP$
admits an extension $\bar f\colon LP\ra LP$.  Assume for a moment  that
$f(x^{(0)}_1)=(x^{(0)}_1)^k\in LP$ for some $k\in \Z$. From the definition 
of $LP$ it follows then that $\bar f$ can be defined by setting $\bar f(r)=kr$ for all 
$r\in \Q$. Next,  notice  that since $P$ is perfect $f(P)\subseteq [LP, LP]=P$. 
Combining these observations we get that the existence of $\bar f$ follows from

 \begin{lemma}
\label{P MORPHISMS}
If $g\colon P\rightarrow P$ is any homomorphism then 
$g(x^{(0)}_1)=(x^{(0)}_1)^k$ for some $k\in \Z$
\end{lemma}

Recall the group $R=P/\langle x_1^{(0)}\rangle$ defined in the proof of 
Lemma \ref{LP TORSION}. Lemma \ref{P MORPHISMS} 
will follow if  we show that for any homomorphism $g\colon P\rightarrow R$
the element $x^{(0)}_1$ is in the kernel of $g$.  
In the proof of  Propositon \ref{P PERFECT} we also defined the group
$$P_1=\langle x^{(1)}_1, x^{(1)}_2 \   | \   [x^{(0)}_1, x^{(1)}_i]=1, \ i=1,2\rangle$$
Since $x^{(0)}_1$ is not in the kernel of the map $P_1\rightarrow P$
it is enough to show that for any $g\colon P_1\rightarrow R$
we have $x^{(0)}_1\in \ker g$. Furthermore, since $R=\ralim_n R_n$ (see \ref{LP TORSION}), 
and $P_1$ is a finitely presented group it  suffices to prove that $g(x^{(0)}_1)=1$
for all  $g\in \Hom(P_1, R_n)$. Finally, notice that by the definition 
of $P_1$ the elements $x^{(1)}_1$, $x^{(1)}_2$ commute with 
their commutator $x^{(0)}_1=[x^{(1)}_1, x^{(1)}_2]$. These observations  
and the presentation of  $R_n$ show that Lemma \ref{P MORPHISMS} 
is a special case of

\begin{lemma}
\label{COMMUTE}
Let  $F_1$, $F_2$ be two free groups,  and let $u_i$ be a word in $F_i$
which is not a proper power. Let  $G$ be the quotient group of $F_1\ast F_2$
by the normal subgroup generated by $[u_1, u_2]$. If $x, y\in G$  are elements 
commuting with $[x, y]$ then $[x, y]=1$.    
\end{lemma}  

\begin{Proof} 
Consider the map $ h\colon G\rightarrow F_1\oplus F_2$. Its kernel $K$  is a free group
 whose set of generators can be described as follows. Let $S_i$ be a set of 
representatives of  cosets of $\langle u_i\rangle\backslash F_i$. Then the generators of 
$K$ are all commutators $[v_1,v_2]$ where $v_1\in S_1$ represents a coset other than 
$\langle u_1\rangle$ and $v_2$ is any nontrivial element of $F_2$, or  $v_2\in S_2$ 
represents a coset  different from $\langle u_2\rangle$, and $v_1$ is a nontrivial 
element  $F_1$. To see this recall \cite[Prop. 4, p. 6]{MR82c:20083} that the kernel $K'$ 
of the map $F_1\ast F_2\ra F_1\oplus F_2$ is a free group whose generators are all 
commutators  $[v_1, v_2]$ where $v_i\in F_i$  and $v_i\neq 1$. The group $K$ is obtained
as the quotient of $K'$ by its normal subgroup generated by the set 
$\{ w^{-1}[u_1, u_2]w \  | \  w\in F_1\ast F_2 \}$. This is equivalent to imposing the following 
relations in $K'$:
$$[u_1w_1, u_2w_2]=[u_1w_1, w_2][w_2, w_1][w_1, u_2w_2]$$
where $w_i$ is an arbitrary element of $F_i$.   The above description of $K$ can be 
derived from here. 
Notice, that  using the above relations any commutator  $[w_1, w_2]$ such that 
$w_i\in F_i$  can be expressed in terms of generators of $K$ using the formula 
\begin{equation}\label{COM EQ}
[w_1, w_2]=[w_1, s_2][s_2, s_1][s_1, w_2]
\end{equation}
where $s_i\in S_i$ represents the coset $\langle u_i \rangle w_i$ for $i=1,2$. 

Next, take  elements $x, y\in G$ as in the statement of the lemma. Notice that $[x,y]\in K$ 
since otherwise  $h(x), h(y)$ would have to commute with a nontrivial element $h([x, y])=[h(x), h(y)]$ 
in $F_1\oplus F_2$ which is impossible. Furthermore, since centralizers of all nontrivial elements 
in a free group are abelian  and since $x, y$ are in the centralizer of $[x, y]$, if $x, y\in K$ then 
we get $xy=yx$ and the statement of the lemma holds. Therefore we can assume that 
$x\not\in K$ and $[x, y]\in K$.  In this case we can uniquely represent $x$ and $[x,y]$ in the form
$$x=x_1x_2\prod_{i=1}^{p} [a_1, b_1]^{\delta_i}\   \text{ and} \  \   
[x, y]=\prod_{j=1}^{q}[c_j, d_j]^{\delta_j}$$ 
where $x_i\in F_i$,  $x_1x_2\neq 1$,  $[a_i, b_i], [c_j, d_j]$ are generators of $K$, and 
$\delta_i, \delta_j= \pm 1$.   
We can also assume that $[x, y]$ is represented by a cyclically reduced word in the free 
group $K$, that is  $[c_1, d_1]^{\delta_1}\neq [c_q, d_q]^{-\delta_q}$. 
Consider the element
\begin{equation}\label{CONJ EQ}
x_2^{-n}x_1^{-n}[x,y]x_1^n x_2^n=\prod_j([x_2^n, c_jx_1^n ][c_jx_1^n, d_jx_2^n]
[d_jx_2^n, x_1^n][x_1^n, x_2^n])^{\delta_j}
\end{equation}
Commutativity of $x$  and $[x, y]$ implies that  for any $n\in \Z$ this element  
is conjugate to $[x, y]$ in $K$. 
We will show that this is not possible unless $[x,y]=1$. One can check that the following holds.
\begin{lemma}\label{CONJ COM}
Let $F$ be a free group and let $w, v$ be words in $F$. If $w$ is cyclically reduced and $v$ is 
conjugated in $F$ to $w$ then all generators of $F$ appearing in $w$ must appear in $v$.
\end{lemma} 
\noindent We apply it to $w=[x,y]$ and $v=x_2^{-n}x_1^{-n}[x,y]x_1^n x_2^n$. The commutators
appearing on the right hand side of  formula 2.8 are not generators of the group K. 
Each of them, however,  can be written as a product of generators of $K$ using formula 2.7. 
By lemma \ref{CONJ COM} all commutators $[c_j, d_j]$ must appear among these 
generators. One can check however that (since $u_1, u_2$ are not proper powers) if $x_1\neq 1$,  $x_2\neq 1$,  and $n$
is large enough this can happen only if $x_1=c^{-1}u_1^kc$, $x_2=d^{-1}u_2^ld$ and $[x,y]=[c,d]^m$ for some $c\in F_1$,  $d\in F_2$, and $k, l, m\in \Z$.  By inspection, in this case $x_1^{-n}x_2^{-n}[x,y]x_1^n x_2^n$ is not conjugated to $[x,y]$ unless $m=0$,  and 
$[x,y]=1$. Assume in turn that e.g. $x_2=1$. Then we have 
$$x_1^{-n}[x,y]x_1^n=\prod_{j=1}^q([c_jx_1^n, d_j][d_j, x_1^n])^{\delta_j}$$
Again, combining this with formula 2.7 we get an  expression of $x_1^{-n}[x,y]x_1^n$
as a product of generators of  $K$. In order  for $[c_j, d_j]$ to appear among these generators 
for large $n$ we must have $x_1=c^{-1}u_1^kc$ and $[x, y]=\prod_j[c, d_j]^{\delta_j}$  for 
some $c\in F_1$, $k\in \Z$. As before, by inspection we obtain that also in this case 
$x_1^{-n}[x,y]x_1^n$ cannot be conjugate to $[x, y]$ if $[x,y]\neq 1$.

\end{Proof}


\end{document}